\documentclass[a4paper,leqno,12pt]{amsart}

\usepackage[utf8]{inputenc}
\usepackage{amsmath}
\usepackage{amsmath}\usepackage{amscd}\usepackage{amssymb}\usepackage{mathrsfs}
\usepackage{stmaryrd}\usepackage{pst-node}\usepackage{tikz-cd}
\usepackage[left=3.2cm,right=3.2cm,bottom=3cm,top=3cm]{geometry}

\newtheorem{theorem}{Theorem}[section]

\newtheorem{proposition}{Proposition}[section]
\newtheorem{lemma}{Lemma}[section]
\newtheorem{corollary}{Corollary}[section]

\newtheorem{remark}{Remark}[section]
\newtheorem{conjecture}{Conjecture}[section]
\newtheorem{question}{Question}[section]

\newcommand{\Ric}{\mathrm{Ric}}

\newcommand{\del}{\nabla}
\newcommand{\laplace}{\triangle}
\def\endproof{{\hfill $\square$}\medskip}

\AtEndDocument{\bigskip{\footnotesize
		\textsc{Macao Center for Mathematical Sciences, Macau University of Science and Technology, Macao, China}\par
		\textit{E-mail address: }\texttt{wxmai@must.edu.mo}\par
		\medskip
		\textsc{School of Mathematical Sciences, Xiamen University, Xiamen, 361005, China}\par
		\textit{E-mail address: }\texttt{oujianyu@xmu.edu.cn}
		
}}

\begin{document}
\title[Liouville Theorem on CSC shrinker]{Liouville Theorem on Ricci shrinkers with constant scalar curvature and its application}
\author{Weixiong Mai, Jianyu Ou}

\maketitle

\begin{abstract}
In this paper we consider harmonic functions on gradient shrinking Ricci solitons with constant scalar curvature. {A Liouville theorem is proved without using gradient estimate :} any bounded harmonic function is constant on gradient shrinking Ricci solitons with constant scalar curvature. As an application, we show that the space of harmonic functions with polynomial growth has finite dimension.
\end{abstract}

\let\thefootnote\relax
\footnotetext{2020 Mathematics Subject Classification 53C25, 53C21}

\footnotetext{Keywords: Liouville theorem, Ricci soliton, frequency function, constant scalar curvature, harmonic function with polynomial growth}

\section{Introduction}

The classical Liouville theorem which is named after Joseph Liouville states that any bounded holomorphic function on the complex plane $\mathbb C$ must be constant, which also holds for harmonic functions on $\mathbb R^n.$
It has a huge impact crossing many fields such as complex analysis, partial differential equations, geometry and so on. 
In 1975, S. T. Yau {\cite{Yau1}} proved that any bounded harmonic function must be constant on an $n$-dimension manifold with nonnegative Ricci curvature, and later, he and S. Y. Cheng {\cite{ChengYau}} showed a gradient estimate on the manifold with Ricci curvature lower bounded to give an effective version of the Liouville theorem, where the gradient estimate method plays a crucial role in the study of harmonic function on complete manifolds. Since Yau's seminal paper, the Liouville type property and the gradient method have been generalized to many partial differential equations on manifolds (see e.g. \cite{Li, CM3} and the references therein for more information).

{
Let $(M^n,g)$ be a complete Riemannian manifold. Denote by $\mathcal H_d(M)$ the space of harmonic functions with polynomial growth at most $d$:
\begin{align*}
    \mathcal H_d(M) = \{u:u(x) \text{ is harmonic on } M^n; \sup_{B_p(r)}|u(x)|\leq  C(1+r)^d \text{ for some }C>0 \},
\end{align*}
where $B_p(r)$ is the geodesic ball with radius $r$ centered at $p.$
In {\cite[1987]{Yau2}}, Yau gave a conjecture: 
\begin{conjecture}{\label{Yau}}
Let $M$ be a complete manifold with nonnegative Ricci curvature, then $\text{dim}(\mathcal H_d(M))<\infty$ for any $d>0.$
\end{conjecture}
This conjecture can be regard as a more general Liouville property. 
The conjecture was firstly proved by P. Li and L.-F. Tam for the case $d=1$ in \cite{Li-Tam89}, and for the $2$-dimensional case ($n=2$) in \cite{Li-Tam91}. 
 T. H. Colding and W. P. Minicozzi ({\cite{CM1,CM2}}) completely proved the conjecture in 1997. In {\cite{CM1}}, they proved the result of finite dimension of $\mathcal H_d(M)$ under weaker assumptions: $M$ is a complete manifold satisfying a volume doubling bound and a scale-invariant Poincare inequality. One can see more in {\cite{CM3}}. In fact, the study of $\mathcal H_d(M)$ has been further developed, see e.g. \cite{Xu,Ding,Huang} for the existence of non-constant functions in $\mathcal H_d(M)$, and see \cite{Liu, Ni} for the study of Yau's conjecture in complete K\"ahler manifolds.
}


{In this paper, we consider a Liouville theorem and the space of harmonic function with polynomial growth on complete gradient Ricci solitons with constant scalar curvature. In fact, the study of Liouville type theorem and {\bf Conjecture \ref{Yau}} in gradient shrinking Ricci solitons have been also considered by some researchers, recently. To state the results, we introduce some notations.}

A complete Riemannian manifold $(M^n, g)$ is a gradient shrinking Ricci soliton if there exists a smooth function $f$ on $M$ {satisfying} the equation
\begin{eqnarray*}
\Ric+\nabla\nabla f={\frac{1}{2}}g,
\end{eqnarray*}
where $\Ric$ is the Ricci tensor and $\nabla\nabla f$ is the Hessian of the function $f$.
The function $f$ is called a potential function of the gradient shrinking Ricci soliton. Ricci soliton can be seen as natural extensions of Einstein metrics. Obviously, if $f$ is a constant function, the gradient Ricci soliton is simply an Einstein manifold. On the other hand, gradient Ricci solitons are also self-similar solutions to Hamilton's Ricci flow and play an important role in the study of formation of singularities in the Ricci flow. We {refer to} {\cite{Cao}} for a nice survey on the subject.

{For shrinking Ricci solitons, many people studied the $f$-harmonic functions, that is the function $u$ such that $\laplace u-<\nabla u, \nabla f>=0$. H. Ge and  S. Zhang {\cite{GZ}} showed that any positive $f$-harmonic function on shrinkers is constant. However, for the standard harmonic function on gradient shrinking Ricci solitons, we do not know whether any bounded harmonic function is constant. In {\cite{MS}} O. Munteanu and N. Sesum proved the harmonic function $u$ with $\int_{M}|\nabla u|^2dv<\infty$ on gradient shrinking Kahler-Ricci soliton is constant.  }

{As an analogue of Yau's conjecture ({\bf Conjecture \ref{Yau}}), an interesting question is following:
\begin{question}\label{Yau_s}
Let $M$ be a complete shrinking Ricci soliton, then $\text{dim}(\mathcal H_d(M))<\infty$ for any $d>0.$
\end{question}
}
{For this question, the following results are known. O. Munteanu and J. Wang {\cite{MW}} showed that the space of holomorphic functions with polynomial growth on Kahler-Ricci shrinkers has finite dimension. Recently, J.-Y. Wu and P. Wu {\cite{JWPW}} proved that the space of harmonic functions with polynomial growth on complete non-compact shrinkers with $R(x)d^2(x,o)\leq c_0$ has finite dimension, where $d(x,o)$ is the distance function from point $x\in M$ to a fixed poiont $o\in M$.}

{ In this paper we will prove a Liouville theorem and affirmatively answer {\bf Question \ref{Yau_s}} on a gradient shrinking Ricci soliton with constant scalar curvature. Our main results are stated as follows.
\begin{theorem}{\label{Liouville}}
Let $(M, g, f)$ be a complete non-compact gradient shrinking Ricci soliton with constant scalar curvature. Any bounded harmonic function on $M$ is constant.
\end{theorem}

Since K\"ahler-Ricci solitons are also Ricci solitons, we have the following corollary.
\begin{corollary}
All bounded harmonic (holomophic) functions on complete gradient shrinking K\"ahler-Ricci solitons with constant scalar curvature are constant.
\end{corollary}

Consequently, using Theorem {\ref{Liouville}}, we obtain a uniform doubling property of harmonic functions with polynomial growth. By the {argument} in {\cite{MW}}, we {prove} the following theorem.
\begin{theorem}{\label{Dimension}}
Let $(M, g, f)$ be a complete non-compact gradient shrinking Ricci soliton with constant scalar curvature, then $\text{dim}(\mathcal H_d(M))<\infty$ for any $d>0.$
\end{theorem}
}

Recently, Colding and Minicozzi {\cite{CM4}} considered the solution of the heat equation which is often called caloric functions. We say $u\in \mathcal{P}_d(M)$ if $u$ is an ancient caloric function and for some $p\in M$ and a constant $C_n$ such that
$$\sup_{B_p(r)\times[-r^2,0]}|u|\leq C_u(1+r)^d,\ \forall r>0.$$
They proved that if a complete manifold $M$ has polynomial volume growth and $k$ is a nonnegative integer, then
$$\dim \mathcal{P}_{2k}(M)\leq\sum^k_{i=0}\dim \mathcal{H}_{2i}(M).$$
By Theorem {\ref{Dimension}}, we obtain the following corollary.

\begin{corollary}
Let $(M, g, f)$ be a complete non-compact gradient shrinking Ricci soliton with constant scalar curvature, then $\text{dim}(\mathcal P_d(M))<\infty$ for any $d>0.$
\end{corollary}

\begin{remark}
{
In many studies of Yau's conjecture ({\bf Conjecture \ref{Yau}}), the method of frequency function plays a crucial role (see \cite{CM1,CM2,Xu}). In this paper our argument, including the proofs of Theorem \ref{Liouville} and Theorem \ref{Dimension}, is also based on the method of frequency function. However, instead of using the frequency function given in the above references, our frequency function is from \cite{Zhu} (see also \cite{Ou}), which has its own advantage in computation.
It should be also noted that J. Ou in {\cite{Ou}} proved a monotonicity result for the used frequency function on gradient shrinking Ricci solitons with constant scalar curvature, which plays a role in our proof of Theorem \ref{Dimension}.}
\end{remark}

\begin{remark}
In fact, a lot of attention has been paid to the study of gradient Ricci solitons with constant scalar curvature in recent years.

Peterson and Wylie {\cite{PW1}} defined a gradient Ricci soliton $(M,g,f)$ to be rigid if it is a flat bundle $N\times_{\Gamma}\mathbb{R}^k$, where $N$ is an Einstein manifold, $\Gamma$ acts freely on $N$ and by orthogonal transformations on $\mathbb{R}^k$, and $f=\frac{\lambda}{2}|x|^2$ on $\mathbb{R}^k$. They \cite{PW1} proved that a gradient Ricci soliton is rigid if and only if it has constant scalar curvature and is radially flat, i.e. $Rm(\del f, \cdot, \del f, \cdot)=0$. 

 H.-D. Cao proposed the following conjecture:

\begin{conjecture}[H.-D. Cao] \label{Caoconjecture}
Gradient shrinking Ricci solitons with constant scalar curvature are rigid.
\end{conjecture}

Fern\'andez-L\'opez and Garc\'ia-R\'io \cite{FG} investigated gradient shrinking Ricci solitons with constant scalar curvature using isoparametric theory, they proved that the constant scalar curvature has to be $0,\ 2\lambda,\ ...,\ (n-1)\lambda$, or $n\lambda$. They \cite{FG} classified four-dimensional and six-dimensional gradient K\"ahler-Ricci solitons of constant scalar curvature. J.-Y. Wu, P. Wu and Wylie  \cite{WWW} classified four-dimensional gradient K\"ahler-Ricci solitons with constant scalar curvature, independently. Finally, X. Cheng and D. Zhou {\cite{ChengZhou}} finished the classification of 4-d gradient Ricci shrinking solitons with constant scalar curvature.

{
Whether {\bf Conjecture {\ref{Caoconjecture}}} for high dimension holds or not is still unknown. 
In some sense Theorem \ref{Liouville} and Theorem \ref{Dimension} provide some positive information on this conjecture since our results are expected if the conjecture has an affirmative answer.}
\end{remark}

{The paper is organized as follows.} We recall some properties of gradient shrinking Ricci solitons in Section 2. We will prove Theorem {\ref{Liouville}} in Section 3. Theorem {\ref{Dimension}} will be proved in Section 4.

{\bf Acknowledgement} The second author thanks Professor Huai-Dong Cao and Professor Fei He for helpful conversations. The first author was supported in part by NSFC Grant No. 11901594. The second author was supported in part by the Fundamental Research Funds for the Central Universities (No. 20720220042).

\section{Preliminaries}

We recall some basic results for gradient shrinking Ricci solitons with constant scalar curvature.
\begin{lemma}{(Hamilton \cite{Ham1})\label{Ham1}}
Let $(M^n,g_{ij},f)$ be a complete gradient shrinking Ricci soliton. Then, we have
\begin{eqnarray*}
\del_i R=2R_{ij}\del_j f
\end{eqnarray*}
and 
\begin{eqnarray*}
R+|\del f|^2-f=C_0
\end{eqnarray*}
for some constant $C_0$. 
\end{lemma}

\begin{lemma} [B.-L. Chen \cite{Chen}, Pigola-Rimoldi-Setti \cite{PRS}]
Let $(M^n, g, f)$ be a complete gradient shrinking Ricci soliton. Then the scalar curvature is nonnegative. Moreover, the scalar curvature is positive unless $(M^n, g, f)$ is the Gaussian soliton $(\mathbb{R}^n, g_0, \frac{\lambda}{2}|x|^2)$.
\end{lemma}

In constant scalar curvature case, {\cite{FG}}, Fern\'andez-L\'opez and Garc\'ia-R\'io gave the following proposition:
\begin{proposition}[Fern\'andez-L\'opez\& Garc\'ia-R\'io \cite{FG} \label{R1}]
Let $(M, g, f)$ be a gradient shrinking Ricci soliton with constant scalar curvature. Then the scalar curvature $R\in \{0, {1\over2}, 1, ..., {n\over 2}\}$. Moreover, no complete gradient shrinking Ricci soliton may exist with $R={1\over 2}$.
\end{proposition}

\begin{remark}
If $R={n\over2}$, we obtain $M$ is n-dimension Einstein manifold. In fact when $(M, g, f)$ is a complete non-compact gradient shrinking Ricci soliton with constant scalar curvature $R$, then $R\in\{0, 1, {3\over2}, 2,..., {{n-1}\over 2}\}$.
\end{remark}

H.-D. Cao and Zhou \cite{CaoZhou} proved that the potential function of gradient shrinking Ricci solitons grows quadratically,
\begin{lemma}[Cao-Zhou \cite{CaoZhou}]\label{cao-zhou}
Let $(M, g, f)$ be a complete non-compact gradient shrinking Ricci soliton. Then there exist constants $c_1,c_2$ such that
\begin{equation*}
\begin{split}
\frac{\lambda}{2}(r(x)-c_1)^2\leq f(x)\leq\frac{\lambda}{2}(r(x)+c_2)^2
\end{split}
\end{equation*}
where $r(x)=d(x_0,x)$ is the distance function from some fixed point $x_0\in M$.
\end{lemma}
{
Note that if we normalize $f$ to $f_0$ by adding  the constant $C_0$, then we have
\begin{eqnarray}{\label{prop1}}
R+|\del f_0|^2=f_0
\end{eqnarray}
Moreover, if the scalar curvature $R$ is a constant and we set the new $f=f_0-R$, we obtain
\begin{eqnarray}{\label{Rij}}
2R_{ij}\del_j f=\del_i R=0
\end{eqnarray}
and 
\begin{eqnarray}{\label{newf}}
|\del f|^2=f.
\end{eqnarray}
Here $f$ is also a potential function of $M$. 
Let $\rho=2\sqrt{f}$, then $|\del \rho|\equiv 1$. {In the rest of the paper, $f$ is the one satisfying (\ref{newf}) and (\ref{Rij}).}
}

Let $(M, g,f)$ be a shrinking Ricci soliton. Set $D_t:=\{x\in M| f(x)\leq t\}$.
In {\cite{ChengZhou}}, Cheng and Zhou have already obtained the volume growth of $D_t$. Here we prove it again:

\begin{theorem}[Cheng-Zhou {\cite{ChengZhou}\label{Volumegrowth}}]
Let $(M^n, g, f)$ be a complete shrinking Ricci soliton with constant scalar curvature $R$ as above. $\Omega(r):=\{x\in M|\rho(x)\leq r\}.$ Then $Vol(\Omega(r))=cr^{n-2R}$ for some constant $c$.
\end{theorem}
{\bf{Proof:}} 
By the method in {\cite{CaoZhou}}:

Let $V(r)=\int_{\Omega(r)}dV$. By co-area formula, $V(r)=\int^r_0\int_{\partial \Omega(s)}{1\over |\del \rho|}dAds$. That means
$$V'(r)=\int_{\partial \Omega(r)}{1\over |\del \rho|}dA=\int_{\partial \Omega(r)}dA.$$

First we have
\begin{eqnarray*}
(n-2R)V(r)&=&nV(r)-2\int_{\Omega(r)}RdV\\
&=&2\int_{\Omega(r)}\laplace fdV\\
&=&2\int_{\partial \Omega(r)}\del f\cdot{\del \rho\over |\del \rho|}dA\\
&=&2\int_{\partial \Omega(r)}|\del f|dA\\
&=&rV'(r).
\end{eqnarray*}
That means ${V(r)\over V'(r)}={n-2R\over r}$, i.e. $\log V(r)=(n-2R)\log r+\log(V(1)).$ So $V(r)=cr^{n-2R}$, where $c=V(1)$ is a constant.
\endproof

\begin{remark}
By Theorem {\ref{Volumegrowth}}, we obtain $Vol(D_t)=ct^{{n\over2}-R}$ for some positive constant $c$. The area of $\partial D_t$ is $Area(D_t)=c({n-2R}) t^{{n\over2}-R-{1\over2}}$.
\end{remark}

Now let us introduce some notations.
For a harmonic function $u$, we define $$H(t)=\int_{D_t}u^2(t-f)^{\alpha}dv,$$ and
$$J(t)=\int_{D_t}|\del u|^2(t-f)^{\alpha+1}$$
with $\alpha\geq 0.$ In particular, we denote by
$$h(t)=\int_{D_t}u^2dv$$
for $H(t)$ with $\alpha=0.$
Then we define the frequency function as $$N(t)=\frac{J(t)}{H(t)}.$$
When $\alpha\geq 2,$ a direct computation yields
\begin{align}\label{HJ}
H^\prime(t) =\frac{\alpha+\frac{n}{2}}{t}H(t) +\frac{\alpha}{t}\int_{D_t} R(t-f)^{\alpha-1} u^2dv + \frac{2}{t(\alpha+1)} J(t) -\frac{1}{t}\int_{D_t} R(t-f)^\alpha u^2dv
\end{align}
and one can show that
$$H(t)\leq t^{\alpha}h(t)$$ and for any $0<t<s$,
$$h(t)\leq {H(s)\over{(s-t)^{\alpha}}}.$$ Moreover, in \cite{Ou} it is also shown that $t^{\sqrt n-1}N(t)$ with $\alpha\geq2$ is nondescreasing for gradient shrinking Ricci solitons with constant scalar curvature, which plays a role in the present paper.
We refer to \cite{Ou} for general discussions on $N(t).$

{In this paper we will focus on the case of} constant scalar curvature $R$. {As mentioned previously,} we choose the potential function $f$ with $|\nabla f|^2=f$.
By (\ref{HJ}), we have
\begin{eqnarray}
H'(t)={\alpha+{n\over2}-R\over t}H(t)+{2\over{(\alpha+1)t}}J(t),
\end{eqnarray}
and 
\begin{eqnarray}{\label{diedai}}
{d\over{dt}}(\ln H(t)))={{\alpha+{n\over2}-R}\over t}+{2\over(\alpha+1)t}N(t).
\end{eqnarray}

\begin{proposition}{\label{P1}}
If $N(t)\geq m>0$, then ${d\over dt}(t^{-{2m\over{\alpha+1}}}H(t)){>}0$.
\end{proposition}

{\bf{Proof:}}
By (\ref{diedai}), we see that ${d\over{dt}}(\ln H(t)))> {2\over(\alpha+1)t}N(t)\geq{2m\over(\alpha+1)t}$. Thus we finish the proof.
\endproof

A harmonic function $u$ is said to be with polynomial growth of order $2d>0$ if
\begin{equation}
    \sup_{x\in D_t}|u(x)| \leq C(t+1)^d.
\end{equation}
Note that the above definition coincides with the standard one (cf. $H_{2d}(M)$) by using Cao-Zhou's result {(see Lemma \ref{cao-zhou}).}

\section{Liouville theorem on Ricci shrinkers with constant scalar curvature}

By the definition of $N(t)$ with $\alpha=0,$ one has
\begin{align}\label{N0}
N(t) =\frac{J(t)}{h(t)}= \frac{\int_{D_t}|\del u|^2(t-f)dv}{\int_{D_t}u^2dv}.
\end{align}

\begin{lemma}\label{N_lim}
For $\alpha=0,$ and $u$ is a non-constant harmonic function satisfying $3\delta >u>\delta,$ where $\delta$ is a positive constant, {there holds}
\begin{align*}
    \lim_{t\to 0} N(t)=0
\end{align*}
and 
\begin{align*}
    \lim_{t\to\infty} N(t)=0.
\end{align*}
\end{lemma}
\noindent{\bf{Proof:}}

For $\alpha=0,$ by a direct computation we have 
\begin{align}\label{N_lim1}
\begin{split}
    N(t) &= \frac{\int_{D_t}|\del u|^2(t-f)dv}{\int_{D_t}u^2dv}\\
    & = \frac{1}{2}\frac{\int_{D_t}\langle \del u^2,\del f \rangle dv}{\int_{D_t}u^2dv}\\
    & = \frac{1}{2}\frac{-\int_{D_t} u^2 \laplace f  dv +\int_{\partial D_t}u^2\langle \del f,\frac{\del f}{|\del f|}\rangle d\sigma }{\int_{D_t}u^2dv}\\
    & = -\frac{1}{2}(\frac{n}{2}-R) + \frac{1}{2}\frac{t\int_{\partial D_t}u^2t^{-\frac{1}{2}}d\sigma}{\int_{D_t}u^2dv}\\
    & = -\frac{1}{2}(\frac{n}{2}-R) + \frac{1}{2}\frac{th^\prime(t)}{h(t)}\\
    & = -\frac{1}{2}(\frac{n}{2}-R) + \frac{1}{2}\frac{(\ln h(t))^\prime}{(\ln t)^\prime}.
    \end{split}
\end{align}
We first consider $\lim_{t\to 0}N(t)=0,$ it suffices to prove
\begin{align*}
    \lim_{t\to 0} \frac{(\ln h(t))^\prime}{(\ln t)^\prime} = \frac{n}{2}-R.
\end{align*}
By LHoptial's rule, we have 
\begin{align}\label{lhopital}
    \lim_{t\to 0}\frac{\ln h(t)}{\ln t^{\frac{n}{2}-R}} = \lim_{t\to 0} \frac{(\ln h(t))^\prime}{(\ln t^{\frac{n}{2}-R})^\prime},
\end{align}
Note that we also have
\begin{align*}
    \delta^2 ct^{\frac{n}{2}-R}\leq \int_{D_t}u^2dv\leq 9\delta^2 ct^{\frac{n}{2}-R},
\end{align*}
which gives 
\begin{align*}
    \ln (\delta^2c) + \ln t^{\frac{n}{2}-R} \leq \ln\int_{D_t}u^2dv\leq \ln ( 9\delta^2c) + \ln t^{\frac{n}{2}-R}.
\end{align*}
For small $t$ we just assume $\ln \int_{D_t}u^2dv<0.$ Therefore,
\begin{align*}
      -\ln (9\delta^2c) + (- \ln t^{\frac{n}{2}-R}) \leq -\ln\int_{D_t}u^2dv\leq -\ln (\delta^2c) + (-\ln t^{\frac{n}{2}-R}),
\end{align*}
which gives 
\begin{align*}
    \frac{-\ln (9\delta^2c)}{(- \ln t^{\frac{n}{2}-R})}+1 \leq \frac{\ln \int_{D_t}u^2}{\ln t^{\frac{n}{2}-R}}\leq \frac{-\ln (\delta^2c)}{(-\ln t^{\frac{n}{2}-R})} +1.
\end{align*}
This leads to $\lim_{t\to 0}\frac{\ln h(t)}{\ln t^{\frac{n}{2}-R}}=1.$ Combining with Lemma \ref{N_lim} and (\ref{lhopital}), we have $\lim_{t\to 0}N(t)=0.$

Similarly, for large $t,$ we have
\begin{align*}
    \frac{\ln (\delta^2c)}{ \ln t^{\frac{n}{2}-R}}+1 \leq \frac{\ln \int_{D_t}u^2dv}{\ln t^{\frac{n}{2}-R}}\leq \frac{\ln (9\delta^2c)}{\ln t^{\frac{n}{2}-R}} +1,
\end{align*}
which implies $\lim_{t\to\infty}\frac{\ln h(t)}{\ln t^{\frac{n}{2}-R}}=1,$ and consequently, $\lim_{t\to\infty}N(t)=0.$
\endproof

Now, we can prove the Theorem {\ref{Liouville}}.

\begin{proposition}\label{liouville}
Let $U$ be a bounded harmonic function, i.e., 
\begin{align*}
    |U(x)|<{B}<\infty,
\end{align*}
where $B$ is a positive constant. Let $u(x)=\frac{\delta U(x)}{B}+\delta+\delta,$ which means $3\delta >u>\delta,$ where $\delta$ is a (small) positive constant. Then we have $u$ is a constant, and consequently, $U$ is a constant.
\end{proposition}

\noindent {\bf Proof:}

By Proposition {\ref{R1}}, we know {in the complete noncompact case} the scalar curvature must be chosen by $\{0, 1, {3\over2}, ..., {n-1\over 2}\}$. We consider the following two cases:

\noindent{\textbf{Case I:}
$\frac{n-2}{2}\leq R{\leq \frac{n-1}{2}}$\\
There are two different proofs for \textbf{Case I}. The first one is similar to that given in \cite[Theorem 11, page 8]{pigola}.
First we suppose that $u$ is not a constant.
Let $K(t) = \int_{D_t} |\del u|^2dv,$ and then $K^\prime(t) = \int_{\partial D_t}|\del u|^2|\del f|^{-1} d\sigma.$
Using integration by parts and Cauchy-Schwarz's inequality, we have
\begin{align*}
    K(t) = & -\int_{D_t} u\laplace u dv +\int_{\partial D_t} u\langle \del u,\frac{\del f}{|\del f|}\rangle d\sigma\\
    = & \int_{\partial D_t} u\langle \del u, \frac{\del f}{|\del f|}\rangle d\sigma\\
    \leq & t^\frac{1}{4}\left(\int_{\partial D_t}u^2dv\right)^\frac{1}{2}\left(\int_{\partial D_t}|\del u|^2 t^{-\frac{1}{2}}d\sigma \right)^\frac{1}{2}\\
    = &t^\frac{1}{4}\left(\int_{\partial D_t}u^2dv\right)^\frac{1}{2}\left(\int_{\partial D_t}|\del u|^2 |\del f|^{-1}d\sigma \right)^\frac{1}{2},
\end{align*}
which gives
\begin{align*}
    -\left(\frac{1}{K(t)}\right)^\prime =\frac{K^\prime(t)}{K^2(t)} \geq & \frac{1}{t^\frac{1}{2}\int_{\partial D_t}u^2dv}\\
    \geq & \frac{1}{9\delta^2 c({n-2R})t^{\frac{n}{2}-R}},
\end{align*}
where we have used $\int_{\partial D_t}d\sigma =c(n-2R)t^{\frac{n}{2}-R-\frac{1}{2}}$ and $3\delta >u>\delta.$
Integrating from $t_1$ to $t_2$, we have
\begin{align*}
    \frac{1}{K(t_1)}> \frac{1}{K(t_1)}-\frac{1}{K(t_2)}\geq \frac{1}{9c\delta^2} \int_{t_1}^{t_2}\frac{1}{t^{\frac{n}{2}-R}}dt.
\end{align*}
Since $\frac{n-2}{2}\leq R\leq \frac{n-1}{2}$, we have $\frac{1}{2}\leq\frac{n}{2}-R \leq 1,$ and consequently, $\frac{1}{t}\leq \frac{1}{t^{\frac{n}{2}-R}}\leq \frac{1}{t^\frac{1}{2}}.$
Therefore, by letting $t_2\to\infty,$ the right-hand side of the above inequality tends to infinity, which implies that $K(t_1)=0.$ We can conclude that $u$ is a constant, which gives a contradiction.

Alternatively, we can also prove \textbf{Case I} by Lemma \ref{N_lim} immediately. By the definition of $N$ with $\alpha=0$ (see (\ref{N0})), we have
\begin{align*}
    \frac{t}{2}\int_{D_\frac{t}{2}}|\del u|^2dv\leq  \int_{D_t}|\del u|^2dv=J(t) = N(t)h(t) \leq 9\delta^2 Vol(D_1)N(t) t^{\frac{n}{2}-R},
\end{align*}
which implies
\begin{align*}
    \int_{D_\frac{t}{2}}|\del u|^2dv \leq C N(t)t^{\frac{n}{2}-R-1}.
\end{align*}
Since $-\frac{1}{2}\leq \frac{n}{2}-R-1\leq 0$, and $\lim_{t\to\infty}N(t)=0$ by Lemma \ref{N_lim}, we have
\begin{align*}
    \int_M|\del u|^2dv = 0
\end{align*}
by letting $t\to \infty.$
Hence, $u$ is a constant.
\bigskip

\noindent{\textbf{Case II:} $0\leq R\leq {{n-3}\over2}$\\

Assume that $u$ is not a constant.
By (\ref{N_lim1}), for such $u$ we also have
\begin{align*}
    h^\prime(t) = \frac{\frac{n}{2}-R}{t}h(t) +\frac{2}{t}J(t),
\end{align*}
which gives
\begin{align*}
    J(t)=\frac{1}{2}(th^\prime(t)-(\frac{n}{2}-R)h(t)).
\end{align*}
Then
\begin{align*}
    2t^{-(\frac{n}{2}-R+1)}J(t)=t^{-(\frac{n}{2}-R)}h^\prime(t) - (\frac{n}{2}-R)t^{-(\frac{n}{2}-R+1)}h(t) = (t^{-(\frac{n}{2}-R)}h(t))^\prime.
\end{align*}
Since $t^{-(\frac{n}{2}-R)}h(t)$ is monotone increasing and $3\delta>u>\delta,$ we have
\begin{align*}
    c\delta^2\leq t^{-(\frac{n}{2}-R)}\int_{D_t}u^2dv \leq 9c\delta^2,
\end{align*}
which implies 
\begin{align*}
    c\delta^2\leq \lim_{t\to 0}t^{-(\frac{n}{2}-R)}h(t)\leq 9c\delta^2
\end{align*}
and 
\begin{align*}
    c\delta^2\leq \lim_{t\to \infty}t^{-(\frac{n}{2}-R)}h(t)\leq 9c\delta^2
\end{align*}
Note that
\begin{align*}
    t^{-(\frac{n}{2}-R)}h(t)-\epsilon^{-(\frac{n}{2}-R)}h(\epsilon) =\int_\epsilon^t (s^{-(\frac{n}{2}-R)}h(s))^\prime ds = 2\int_{\epsilon}^t s^{-(\frac{n}{2}-R+1)}J(s) ds.
\end{align*}
By letting $\epsilon\to 0$ and $t\to \infty$ we have
\begin{align*}
    0<C_1 \leq \int_0^\infty s^{-(\frac{n}{2}-R+1)}J(s)ds\leq C_2<\infty.
\end{align*}
Note that
\begin{align*}
    \int_0^\infty s^{-(\frac{n}{2}-R+1)}J(s) ds & = \int_0^\infty \int_{D_s}s^{-(\frac{n}{2}-R+1)}|\del u(x)|^2 (s-f(x))dv ds\\
    & = \int_0^\infty \int_{M} s^{-(\frac{n}{2}-R+1)}|\del u(x)|^2\chi_{\{x:f(x)<s\}}(x) (s-f(x))dvds\\
    & = \int_{M}|\del u(x)|^2 \int_0^\infty  s^{-(\frac{n}{2}-R+1)}\chi_{\{s:s>f(x)\}}(s) (s-f(x))dsdv\\
    & = \int_{M}|\del u(x)|^2 \int_f^\infty  s^{-(\frac{n}{2}-R+1)} (s-f(x))dsdv\\
    & = \int_M |\del u(x)|^2 \left(\int_f^\infty s^{-(\frac{n}{2}-R)}ds-f(x)\int_f^\infty s^{-(\frac{n}{2}-R+1)}ds \right)dv\\
    & = \int_M |\del u(x)|^2 \left(\frac{s^{-(\frac{n}{2}-R)+1}}{-(\frac{n}{2}-R)+1}\big|_f^\infty-f(x) \frac{s^{-(\frac{n}{2}-R)}}{-(\frac{n}{2}-R)}\big|_f^\infty \right)dv\\
    & = \int_{M}|\del u(x)|^2 \left(\frac{f^{-(\frac{n}{2}-R-1)}}{\frac{n}{2}-R-1}-\frac{f^{-(\frac{n}{2}-R-1)}}{\frac{n}{2}-R}\right)dv\\
    & = \frac{1}{(\frac{n}{2}-R)(\frac{n}{2}-R-1)}\int_M |\del u|^2 f^{-(\frac{n}{2}-R-1)}dv,
\end{align*}
which is bounded from above and below.

Note that
\begin{align*}
    \int_M |\del u|^2  f^{-(\frac{n}{2}-R)+1}dv & = \lim_{t\to \infty} \int_{ D_t} |\del u|^2  f^{-(\frac{n}{2}-R)+1}dv.
\end{align*}
By integration by parts, we have
\begin{align*}
    \int_{ D_t}|\del u|^2 f^{-(\frac{n}{2}-R)+1}dv  = & -\int_{D_t} u\laplace u  f^{-(\frac{n}{2}-R)+1}dv + (\frac{n}{2}-R+1)\int_{ D_t}u\langle \del u,\del  f\rangle  f^{-(\frac{n}{2}-R)}dv\\
    & + \int_{\partial D_t} t^{-(\frac{n}{2}-R)+1}u\langle \del u,\frac{\del f}{|\del f|} \rangle d\sigma.
\end{align*}

By integration by parts again, we have
\begin{align*}
   & \int_{ D_t}\langle \del u^2,\del  f\rangle  f^{-(\frac{n}{2}-R)}dv \\
    = & -\int_{ D_t}  u^2 (\laplace f) f^{-(\frac{n}{2}-R)}dv + (\frac{n}{2}-R)\int_{ D_t} u^2|\del  f|^2 f^{-(\frac{n}{2}-R)-1}dv\\
    & +\int_{\partial  D_t}t^{-(\frac{n}{2}-R)}u^2\langle \del  f,\frac{\del  f}{|\del f|} \rangle d\sigma\\
    = & \int_{\partial  D_t}t^{-(\frac{n}{2}-R)}u^2\langle \del  f,\frac{\del  f}{|\del f|} \rangle d\sigma,
\end{align*}
where we have used $\laplace f=\frac{n}{2}-R$ and $f=|\del f|^2.$

By a combination of the above equalities, we have
\begin{align*}
    &\int_{ D_t}|\del u|^2 f^{-(\frac{n}{2}-R)+1}dv\\
    = &(\frac{n}{2}-R+1) t^{-(\frac{n}{2}-R)}\int_{\partial  D_t}u^2\langle \del  f,\frac{\del  f}{|\del f|} \rangle d\sigma + \int_{\partial D_t} t^{-(\frac{n}{2}-R)+1}u\langle \del u,\frac{\del f}{|\del f|} \rangle d\sigma\\
    = & (\frac{n}{2}-R+1)\int_{\partial D_t}t^{-(\frac{n}{2}-R)+\frac{1}{2}} u^2d\sigma +t^{-(\frac{n}{2}-R)+1} \int_{D_t} |\del u|^2 dv. 
\end{align*}
That is
\begin{align*}
    &\int_{ D_t}|\del u|^2 f^{-(\frac{n}{2}-R)+1}dv
    = (\frac{n}{2}-R+1) \int_{\partial D_t}t^{-(\frac{n}{2}-R)+\frac{1}{2}} u^2d\sigma +t^{-(\frac{n}{2}-R)+1} \int_{D_t} |\del u|^2 dv. 
\end{align*}
Note that
\begin{align*}
    J(t) = N(t)h(t),
\end{align*}
which implies that
\begin{align*}
    t^{-(\frac{n}{2}-R)+1}\int_{D_{\frac{t}{2}}}|\del u|^2 dv \leq C_3N(t),
\end{align*}
where the right-hand side tends to $0$ as $t\to \infty$ or $t\to 0$ given by Lemma \ref{N_lim}. 
Hence,
\begin{align*}
    \int_M |\del u|^2 f^{-(\frac{n}{2}-R)+1}dv =(\frac{n}{2}-R+1) \lim_{t\to \infty} \int_{\partial D_t}u^2 t^{-(\frac{n}{2}-R)+\frac{1}{2}}d\sigma,
\end{align*}
and
\begin{align*}
   0= \lim_{t\to 0} \int_{D_t}|\del u|^2f^{-(\frac{n}{2}-R)+1}dv =(\frac{n}{2}-R+1) \lim_{t\to 0} \int_{\partial D_t} u^2 t^{-(\frac{n}{2}-R)+\frac{1}{2}} d\sigma \geq C_4>0,
\end{align*}
which gives a contradiction.
Hence, $u$ is a constant.
\endproof
\begin{remark}
In fact, for {\textbf{Case I:} $\frac{n-2}{2}\leq R{\leq \frac{n-1}{2}},$} by the result given in \cite[Theorem 11,page 8]{pigola}, one can have that any positive harmonic function is constant, which is  stronger than the conclusion given in Proposition \ref{liouville}.
\end{remark}

Similarly by Corollary 3.2 in {\cite{MS}}, we can also obtain the following corollary about the ends of shrinkers with constant scalar curvature.

\begin{corollary}
Let $(M, g, f)$ be a complete gradient shrinking Ricci soliton with constant scalar curvature. It has at most one nonparabolic end.
\end{corollary}
{\bf{Proof:}} Let $l$ be the number of nonparabolic ends. By the theory of Li and Tam ({\cite{LT}} Theorem 2.1), $l\leq \dim \mathcal{H}_0(M)$ which finish the proof.
\endproof

Then by the result in ({\cite{MS}} Proposition 2.3), we know all the ends of $M$ are nonparabolic if $R\leq a <{n-2\over 2}$. We obtain the following result.

\begin{corollary}
If $(M, g, f)$ is a complete gradient shrinking Ricci soliton with constant scalar curvature $R\leq {n-3\over2}$, then it is connected at infinity.
\end{corollary}

\begin{remark}
In fact, by the splitting Theorem in {\cite{MW2}} (Theorem 1.7), we know that if there are two ends on a complete gradient shrinking Ricci soliton $(M, g, f)$ with constant scalar curvature. Then $M$ is isometric to $\mathbb{N}^{n-1}\times \mathbb{R}$, where $\mathbb{N}$ is an Einstein manifold with positive scalar curvature.
\end{remark}


\section{Application : dimension estimate}

To prove the Theorem {\ref{Dimension}}, we need to show that the following doubling property for harmonic function with polynomial growth.

\begin{proposition}\label{poly}
For any $d>0,$ if $u$ is harmonic with polynomial growth with order $2d,$ i.e., 
$$
    \sup_{x\in D_t}|u(x)|\leq C(t+1)^d,
$$
then for any $1<t<T<\infty,$ there exists a large $\lambda$ such that $$N(t)\leq {(\alpha+1)\lambda^{\sqrt{n}-1} T^{\sqrt n-1}}(d+\epsilon+\frac{\frac{n}{2}-R+\alpha}{2}),$$
with $\alpha\geq2$, $\epsilon>0.$
Moreover, 
$$
\frac{H(t_2)}{H(t_1)} \leq \frac{t_2^L}{t_1^L},
$$
or equivalently, $H(2t)\leq 2^L H(t)$ for $1<t<T,$
where $L=\alpha+\frac{n}{2}-R+2\lambda^{\sqrt n-1}T^{\sqrt n-1}(d+\epsilon +\frac{\frac{n}{2}-R+\alpha}{2}).$
\end{proposition}
\noindent {\bf Proof:} 
We denote $N_u(t)=N(t)$. If the conclusion does not hold, then for any large $\lambda_i$, there exist $u_i$ and $T>t_i>1$ such that $N_{u_i}(t_i)=N_i(t_i)>(\alpha+1)\lambda_i^{\sqrt n-1}T^{\sqrt n-1}(d+\epsilon+\frac{\frac{n}{2}-R+\alpha}{2})$.
Note that when $\alpha\geq 2,$ $t^{\sqrt{n}-1}N(t)$ is nondescreasing (see \cite{Ou}). Then we have
\begin{eqnarray}
 t^{\sqrt n-1}N_i(t)\geq t_0^{\sqrt n-1}N(t_0)>(\alpha+1)\lambda_i^{\sqrt n-1} T^{\sqrt n-1}(d+\epsilon+\frac{\frac{n}{2}-R+\alpha}{2}),
\end{eqnarray}
which holds for all $t>t_i$. In particular, we  have
\begin{eqnarray*}
N_i(t)> (\alpha +1)(d+\epsilon+\frac{\frac{n}{2}-R+\alpha}{2}),
\end{eqnarray*}
for $\lambda_i T>t>T$.
Then applying Proposition \ref{P1} with $m=(\alpha +1)(d+\epsilon+\frac{\frac{n}{2}-R+\alpha}{2}),$ we have
$t^{-2(d+\epsilon+\frac{\frac{n}{2}-R+\alpha}{2})}H(t)$ is nondescreasing for $\lambda_i T>t>T$.
Consequently, we have
$$
    t^{-2(d+\epsilon+\frac{\frac{n}{2}-R+\alpha}{2})}H_i(t)\geq T^{-2(d+\epsilon+\frac{\frac{n}{2}-R+\alpha}{2})}H_i(T)
$$
for $T<t<\lambda_i T.$

Then, by the volume of $D_t$ for the constant scalar curvature case (i.e., $Vol(D_t)=ct^{\frac{n}{2}-R}$), there holds
\begin{eqnarray*}
ct^{\alpha+\frac{n}{2}-R}\sup_{D_t}u_i^2 &=& t^\alpha Vol(D_t)\sup_{D_t}u_i^2\\
&\geq&  \sup_{D_t}u_i^2\int_{D_t}(t-f)^\alpha \\
&\geq& H_i(t)\\
&\geq & T^{-2(d+\epsilon+\frac{\frac{n}{2}-R+\alpha}{2})}H_i(T)t^{2d+2\epsilon+{\frac{n}{2}-R+\alpha}}\\
&\geq & T^{-2(d+\epsilon+\frac{n-R+\alpha}{2})}H_i(T)t^{2d+2\epsilon+\frac{n}{2}-R+\alpha},
\end{eqnarray*}
which implies 
\begin{eqnarray*}
c\sup_{D_t}u_i^2 \geq T^{-2d-2\epsilon-\frac{n}{2}+R-\alpha}H_i(T) t^{2d+2\epsilon},
\end{eqnarray*}
for $\lambda_i T>t>T$.

Then we have
\begin{align*}
    cC^2t^{2d} &>T^{-2d-2\epsilon-\frac{n}{2}+R-\alpha}t^{2d+2\epsilon}H_i(T)
\end{align*}
for $\lambda_i T>t>T.$
Therefore
\begin{align*}
    \widetilde C ((\lambda_i -1)T)^{-2\epsilon}>\widetilde C t^{-2\epsilon}>H_i(T),
\end{align*}
for $\lambda_i T>t>(\lambda_i-1)T.$
Since $\lim_{i\to \infty}\lambda_i=\infty,$ we have $\lim_{i\to\infty} H_i(T)=0,$ which implies $\lim_{i\to\infty}u_i=0.$ Thus, for sufficiently large $i,$ we have that $u_i$ is uniformly bounded. By Theorem \ref{Liouville}, we have $u_i$ is a constant, and then $N_i(t)\equiv 0,$ which contradicts to $N_i(t_i)>(\alpha+1)\lambda_i^{\sqrt n-1}T^{\sqrt n-1}(d+\epsilon+\frac{\frac{n}{2}-R+\alpha}{2}).$
 Therefore, there exists a $\lambda$ such that
\begin{eqnarray}
N(t)\leq {(\alpha+1)\lambda^{\sqrt n-1}T^{\sqrt n-1}}(d+\epsilon+\frac{\frac{n}{2}-R+\alpha}{2})
\end{eqnarray}
for $1<t<T.$

Consequently, we have
\begin{eqnarray*}
(\ln H(t))^\prime \leq \frac{L}{t},\quad 1<t<T,
\end{eqnarray*}
with $L=\alpha+\frac{n}{2}-R+2\lambda^{\sqrt n-1}T^{\sqrt n-1}(d+\epsilon+\frac{\frac{n}{2}-R+\alpha}{2}).$ Then
$$
\ln\frac{H(t_2)}{H(t_1)} =\int_{t_1}^{t_2} (\ln H(t))^\prime dt\leq \ln\frac{t_2^L}{t_1^L},
$$
which means $\frac{H(t_2)}{H(t_1)} \leq \frac{t_2^L}{t_1^L}$.
\endproof

\bigskip

The proof of Theorem {\ref{Dimension}} is given as follows.

\noindent{\bf{Proof of Theorem {\ref{Dimension}}:}}

 {By Proposition \ref{poly}} the harmonic function $u$ with polynomial growth has the doubling property:
\begin{eqnarray}
\int_{D(5t)}|u|^2dv\leq C(L)\int_{D(t)}|u|^2dv,
\end{eqnarray}
for {$1<t<5t<T$}.

Using the same argument in {\cite{MW}}, we can have the conclusion. For completeness, we include a proof.

By {\cite{CaoZhou}} (see also Lemma \ref{cao-zhou}), $t$ can be chosen such that
$$D(t)\subset B_p(3r)$$
and 
$$B_p(4r)\subset D(5t),$$
where $p$ is a fixed point in $M$ and $r=\sqrt{t}$. Thus, we have
\begin{eqnarray}{\label{double1}}
\int_{B_p(4r_0)}|u|^2dv\leq C_0(L) \int_{B_p(3r_0)}|u|^2dv,
\end{eqnarray}
for some $r_0$ depending only on $L$.

Note that {(\ref{double1})} holds for all harmonic function with polynomial growth on $M$ with some $L$.

Then, by a result of Li {\cite{Li}}, there exists a nontrivial $u_0\in \mathcal H_{d}(M)$ such that
\begin{eqnarray}{\label{double2}}
\int_{B_p(3r_0)}|u_0|^2dv\leq {nV_p(3r_0)\over{\dim{\mathcal H_{d}(M)}}}\sup_{B_p(3r_0)}|u_0|^2,
\end{eqnarray}
where $V_p(r)$ denotes the volume of $B_p(r).$
On the other hand, the Sobolev constant of $B_p(4r_0)$ can be controlled by a constant depending only on $n$ and $r_0$ ({\cite{MW0}}). Then using the Moser iteration for the subharmonic function $|u_0|$, there holds
\begin{eqnarray}{\label{double3}}
\sup_{B_p(3r_0)}|u_0|^2dv\leq {C_1(n,d)\over{V_p(4r_0)}}\int_{B_p(4r_0)}|u_0|^2dv.
\end{eqnarray}

By (\ref{double1}), (\ref{double2}) and (\ref{double3}), we obtain that
\begin{eqnarray*}
\int_{B_p(4r_0)}|u_0|^2dv&\leq& C_0(L) \int_{B_p(3r_0)}|u_0|^2dv\\
&\leq& C_0(L) {nV_p(3r_0)\over{\dim{\mathcal H_{d}(M)}}}\sup_{B_p(3r_0)}|u_0|^2\\
&\leq& C_0(L) {nV_p(3r_0)\over{\dim{\mathcal H_{d}(M)}}}{C_1(n,d)\over{V_p(4r_0)}}\int_{B_p(4r_0)}|u_0|^2dv,
\end{eqnarray*}
which implies $\dim \mathcal H_{d}(M)\leq \widetilde C_0(L)<\infty$.

\end{document}